\documentclass[10pt]{article}          
\usepackage{amssymb} 
\usepackage{amsmath}               
\usepackage{epsfig,float} 
\usepackage{graphicx}
\usepackage{float}
\usepackage{color}
\usepackage{listings}
\usepackage{authblk}
\usepackage[table]{xcolor}
\usepackage{algorithm}
\usepackage{algorithmic}
\usepackage{placeins}
\usepackage{array}
\usepackage{caption}
\usepackage{subcaption}
\usepackage{mathtools}
\usepackage{tikz}
\usetikzlibrary{plotmarks}
\usepackage{pgfplots}
\pgfplotsset{compat=1.14}

\newenvironment{customlegend}[1][]{
	\begingroup
	\csname pgfplots@init@cleared@structures\endcsname
	\pgfplotsset{#1}
}{
	\csname pgfplots@createlegend\endcsname
	\endgroup
}
\def\addlegendimage{\csname pgfplots@addlegendimage\endcsname}
\pgfkeys{/pgfplots/number in legend/.style={
		/pgfplots/legend image code/.code={
			\node at (0.295,-0.0225){#1};
		},
	},
}

\newcolumntype{C}[1]{>{\centering\arraybackslash\hspace{0pt}}b{#1}}

\newenvironment{keywords}
{\begin{trivlist}\item[]{\bfseries Keywords:}\ }
	{\end{trivlist}}

\newenvironment{AMS}
{\begin{trivlist}\item[]{\bfseries AMS:}\ }
	{\end{trivlist}}

\begin{document}                                          
\title{Real-time implementation of an iterative solver for atmospheric tomography}

\author[1]{Bernadett Stadler}
\author[2]{Roberto Biasi}
\author[2]{Mauro Manetti}
\author[1]{Ronny Ramlau}
\affil[1]{Industrial Mathematics Institute, Johannes Kepler University Linz, Altenbergerstra\ss e 69, 4040 Linz, Austria}
\affil[2]{Microgate, Via Waltraud-Gebert-Deeg 3e, 39100 Bolzano, Italy}

\maketitle

\begin{abstract}
	The image quality of the new generation of earthbound Extremely Large Telescopes (ELTs) is heavily influenced by atmospheric turbulences. To compensate these optical distortions a technique called adaptive optics (AO) is used. Many AO systems require the reconstruction of the refractive index fluctuations in the atmosphere, called atmospheric tomography. The standard way of solving this problem is the Matrix Vector Multiplication, i.e., the direct application of a (regularized) generalized inverse of the system operator. However, over the last years the telescope sizes have increased significantly and the computational efficiency become an issue. Promising alternatives are iterative methods such as the Finite Element Wavelet Hybrid Algorithm (FEWHA), which is based on wavelets. Due to its efficient matrix-free representation of the underlying operators, the number of floating point operations and memory usage decreases significantly. In this paper, we focus on performance optimization techniques, such as parallel programming models, for the implementation of this iterative method on CPUs and GPUs. We evaluate the computational performance of our optimized, parallel version of FEWHA for ELT-sized test configurations. 
\end{abstract}


\begin{AMS}
	{65R32, 65Y05, 85-08, 85-10}
\end{AMS}

\section{Introduction}
Due to the irregular mixing of cold and hot air, affected by the sun and wind, the image quality of the new generation of earthbound Extremely Large Telescopes (ELTs) is heavily disturbed. These irregularities lead to an inhomogeneous refractive index of air ,and thus to a distorted wavefront arriving at the telescope pupil. To compensate these rapidly changing optical distortions, the deformations of wavefronts emitted by natural or laser guide stars (NGS or LGS) are measured via wavefront sensors (WFS) and, subsequently, corrected using deformable mirrors (DMs). A DM typically consists of an optical surface to reflect light, deformed in real-time by a set of actuators. This technique is called Adaptive Optics (AO) \cite{Roddier,RoWe96,ElVo09}. Three different AO systems are based on atmospheric tomography, i.e., they require the reconstruction of the turbulence profile in the atmosphere. These AO systems are: Laser Tomography Adaptive Optics (LTAO), Multi Object Adaptive Optics (MOAO) and Multi Conjugated Adaptive Optics (MCAO). LTAO uses a combination of several LGS and NGS with a single DM to sharpen one scientific object of interest. MOAO is based on the same concept, but uses several mirrors that are optimized to sharpen different objects of interest at the same time. In contrast, MCAO uses several DMs, conjugated to different heights, to achieve a high imaging quality in a large field of view.

Mathematically, the underlying atmospheric tomography problem is ill-posed, i.e., there is an unstable relation between measurements and the solution \cite{Davison83,Nat86}. As a consequence, regularization techniques are required. A common way to regularize this problem is the Bayesian framework, because it allows to incorporate statistical information about turbulence and noise. The random variables are typically assumed to be Gaussian, therefore, the maximum a posterior (MAP) estimate is an optimal point estimate for the solution. Detailed information about the systems can be found in \cite{Hammer02,Andersen06,RiElFl00,Puech08,DBB10}. The dimension of the atmospheric tomography problem depends on the number of subapertures of the used WFS and on the number of degrees of freedom of the DMs, which are in general higher for bigger telescopes. Moreover, the solution has to be computed in real-time, leading to a highly non-trivial task for ELT-sized problems. So far, the standard solver for atmospheric tomography is the Matrix Vector Multiplication (MVM), i.e., the direct application of a (regularized) generalized inverse of the system operator. The computational costs of the MVM in hard real-time scale at $\mathcal{O}(n^2)$, where $n$ is the dimension of the AO system. For computing the inverse in soft-real time the computational costs are even higher and scale at $\mathcal{O}(n^4)$. This direct solution method is suitable for small telescopes, however, if $n$ is increasing significantly, as for the ELT, the computational efficiency becomes an issue. In this paper, we focus on a wavelet based iterative method called FEWHA in which the MAP estimate is computed by a preconditioned conjugate gradient method (PCG) \cite{Yu14, HeYu13, YuHeRa13, YuHeRa13b}. There are several further possibilities on how to treat the atmospheric tomography problem, either directly or iteratively, see \cite{Fusco,ElGiVo02,GiElVo02b,GiElVo03,YaVoEl06,GiElVo07,GiEl08,RoCoGrScFu10,ThiTa10,Tallon_et_al_10,RaRo12,RoRa13,RaObRoSa13,SaRa15,RafRaYu16}. The key feature of FEWHA is that the inverse covariance matrix in the frequency domain can be well approximated by a diagonal matrix. Moreover, a dual domain discretization approach is used to obtain a sparse structure for other operators as well. This concept allows an efficient matrix-free representation of all operators, leading to a significant reduction in floating point operations and memory resources. The quality of the method has been studied in \cite{Yu14,HeYu13, YuHeRa13, YuHeRa13b}. In this paper, we focus on the computational performance of FEWHA on different real-time hardware. Real-time implementations for ELT AO systems are also studied in \cite{GPURTC, FPGARTC, MAORY, AOrealtime, GPURTC2, GPURTC3}. Possible real-time architectures have been evaluated within the Greenflash project, see \cite{Greenflash1, Greenflash2}. Based on these investigations we focus on CPUs and GPUs within this work. We utilize common parallel programming models that have been widely studied, e.g., in \cite{Parallel_Computing, Cuda, ProTBB}, to optimize the performance of the method for specific hardware architectures. We demonstrate the performance of our implementations on the test configuration of MAORY, which is an adaptive optics module for the ELT operating in MCAO. The MAORY real-time computer is discussed in \cite{MAORY}.  

The paper is organized as follows: In Section \ref{sec:prelim} we give a short summary of atmospheric tomography and the iterative solver FEWHA. Section \ref{sec:realtimecomp} is devoted to the optimization techniques, which we apply to FEWHA for specific real-time hardware. Section \ref{sec:numerics} contains numerical simulations on certain hardware architectures for the test configuration of MAORY. Additionally, the behavior of the method for varying parameter settings is analyzed.

\section{Preliminaries}\label{sec:prelim}
\subsection{Atmospheric tomography}
Atmospheric tomography is the fundamental problem in many AO systems. Assuming a layered model of the atmosphere, the goal of atmospheric tomography is to reconstruct the turbulent layers from the wavefront sensor measurements \cite{RoWe96}.

The atmospheric tomography problem is defined by
\begin{equation}\label{eq:atmo}
s = (s_g^x, s_g^y)_{g=1}^G = A\phi,
\end{equation}
where $G$ is the number of guide stars, $\phi = (\phi_1, ..., \phi_L)$ denote the $L$ turbulent layers and $s$ the sensor measurements. We consider here the Shack-Hartmann (SH) WFS. The tomography operator $A$ is decomposed into the SH operator $\Gamma$ and a geometric propagation operator $P$ in the direction of the guide star. Hence, for a specific guide star direction $g$ we obtain
\begin{equation*}
s_g = \Gamma_g P_g \phi \quad \text{ for } g=1,...,G.
\end{equation*}

The SH-WFS utilizes an array of little lenses, each focused on a charge-coupled device (CCD) detector plane. The vertical and horizontal shifts of the focal points determine the average slope of the wavefront over the area of the lens, known as subaperture. We define the subaperture grid of $n_s$ subapertures for a telescope with diameter $D$ as a grid of points with equidistant spacing by
\begin{equation*}
\Omega = \bigcup_{0\leq i,j <n_s} \bar{\Omega}_{ij}
\end{equation*}
with
\begin{equation*}
\Omega_{ij} := (x_i,x_{i+1}) \times (x_j,x_{j+1})
\end{equation*}
and
\begin{equation*}
x_i := -D/2 + i  \quad \text{ for } i=1,...,n_s.
\end{equation*}

The SH measurements in a subaperture $\Omega_{ij}$ are modeled by the average slopes of the wavefront aberration $\varphi$ in $\Omega_{ij}$. Assuming that the incoming wavefront is described by a piecewise-continuous bilinear function $\varphi$ with nodal values $\varphi_{ij}$ at the points $\{(x_i ,x_j) : 0 \leq i,j \leq n_s\}$ this reduces to the following equations
\begin{equation*}
s_{ij}^x=\frac{(\varphi_{i,j+1}-\varphi_{i,j})+(\varphi_{i+1, j+1} - \varphi_{i+1,j})}{2},
\end{equation*}
\begin{equation*}
s_{ij}^y=\frac{(\varphi_{i+1,j}-\varphi_{i,j})+(\varphi_{i+1, j+1} - \varphi_{i,j})}{2}.
\end{equation*}
The SH operator $\Gamma$ maps wavefronts to SH-WFS measurements
\begin{equation*}
s = \begin{pmatrix}s^x\\s^y\end{pmatrix}=\begin{pmatrix}\Gamma^x\varphi\\\Gamma^y\varphi\end{pmatrix} \eqqcolon \Gamma \varphi.
\end{equation*}  

The wavefront aberrations $\varphi$ in the direction of a NGS are given by 
\begin{equation*}
\varphi_\theta(x) = (P_\theta^{NGS}\phi)(x):=\sum_{\ell = 1}^L\phi_\ell(x+\theta h_\ell),
\end{equation*}
where $x=(x_1,x_2,0)$ is a point on the aperture, $\theta = (\theta_1,\theta_2,1)$ is the direction of the guide star and $h_\ell$ is the layer height. The wavefront aberrations in the direction of a LGS at a fixed height $H$ are given by  
\begin{equation*}
\varphi_\theta(x) = (P_\theta^{LGS}\phi)(x):=\sum_{\ell = 1}^L\phi_\ell\left((1-\frac{h_\ell}{H})x+\theta h_\ell\right).
\end{equation*}
For the LGS case the so called cone-effect has to be taken into account. For details about the definition of the geometric propagation operator, either for NGS or LGS, we refer to \cite{Fusco}. 

Mathematically, Equation \eqref{eq:atmo} is ill-posed \cite{Davison83}, i.e., there is an unstable relation between the measurements and the solution. Hence, regularization techniques are required. A common procedure in the literature is to formulate the problem in the Bayesian framework, where statistical information regarding the turbulence model and noise can be utilized. Let $\boldsymbol{S}$ and $\boldsymbol{\Phi}$ be random variables that correspond to the vectors of measurements and turbulence layers, respectively. Further, we assume the presence of noise modeled as a random variable $\boldsymbol{\eta}$ and formulate Equation \eqref{eq:atmo} in the Bayesian framework as 
\begin{equation*}
\boldsymbol{S} = A\boldsymbol{\Phi} + \boldsymbol{\eta}.
\end{equation*}

The random variables are typically assumed to be Gaussian, therefore, the maximum a posteriori (MAP) estimate provides an optimal point estimate for the solution, see \cite{Fusco} for details. We obtain  
\begin{equation*}
x_{MAP} = \arg\min_{\phi \in \mathbb{R}^n} \{\Vert\phi\Vert^2_{C_\phi^{-1}} + \Vert s - A\phi\Vert^2_{C_\eta^{-1}}\},
\end{equation*}
where $C_\phi^{-1}$ and $C_\eta^{-1}$ are the inverse covariance matrices of layers $\boldsymbol{\Phi}$ and noise $\boldsymbol{\eta}$, respectively. The norm induced by a symmetric, positive definite matrix $C$ is defined as
\begin{equation*}
\Vert x \Vert_C^2 := (Cx, x).
\end{equation*}
The solution to this minimization problem is given by the solution of the linear system of equations
\begin{equation}\label{eq:MAP}
(A^* C_\eta^{-1}A + C_\phi^{-1})\phi = A^* C_\eta^{-1}s.
\end{equation}
We assume that the layers are zero centered and uncorrelated, which implies a block-diagonal structure of $C_\phi$. The noise covariance matrix $C_\eta$ is given as a block-diagonal matrix with respect to the wavefront sensors.

The size of the atmospheric tomography problem, i.e., the dimension of $A$, is in general higher for larger telescopes. Within the era of new extremely large ground based telescopes solving the atmospheric tomography problem in real-time is a non-trivial task. For that purpose, efficient solution methods are of great interest. The Finite Element Wavelet Hybrid Algorithm (FEWHA), described in the next section, is one possible approach.

\subsection{Finite Element Wavelet Hybrid Algorithm (FEWHA)}
The Finite Element Wavelet Hybrid Algorithm (FEWHA) is an iterative approach to solve Equation (\ref{eq:MAP}). The main idea is to use compactly supported orthonormal wavelets for representing the turbulent layers. The properties of the wavelet decomposition in the frequency domain allow a completely diagonal approximation of the penalty term  $C_\phi$. The atmospheric tomography operator $A$ has a more efficient representation in a finite element domain, where continuous piecewise bilinear functions are utilized to represent the incoming wavefronts and layers. For details about this dual domain discretization approach we refer to \cite{YuHeRa13b}.

By combining these two representations we obtain the following dual domain discretization of the MAP estimate \eqref{eq:MAP} 
\begin{equation}\label{eq:fewha}
(\boldsymbol{W}^{-T}\hat{A}^TC_\eta^{-1}\hat{A}\boldsymbol{W}^{-1} + \alpha D)c=\boldsymbol{W}^{-T}\hat{A}^TC_\eta^{-1}s,
\end{equation}
where $\hat{A}$ is the atmospheric tomography operator in the finite element domain and $\boldsymbol{W}$ is the discrete wavelet transform, which acts as a linear mapping between the finite element and the wavelet domain. The operator $C_\eta^{-1}$ denotes the inverse covariance matrix of the noise and $D$ is a diagonal approximation of $C_\phi^{-1}$ in the frequency domain. As we have an approximation of $C_\phi^{-1}$ we introduce a scalar factor $\alpha$ for tuning the balance between the fitting and the regularizing terms, see \cite{HeYu13}. The vector $c$ is a concatenation of all wavelet coefficients of all turbulence layers and the vector $s$ is the concatenation of all SH sensor measurements from all guide star directions. Equation \eqref{eq:fewha} is a discretization of \eqref{eq:MAP}, i.e., the operators are now matrices.

For the sake of simplicity, we define the left-hand side operator of Equation \eqref{eq:fewha} by
\begin{equation*}
M := (\boldsymbol{W}^{-T}\hat{A}^T\hat{C_\eta}^{-1}\hat{A}\boldsymbol{W}^{-1} + \alpha D)
\end{equation*}
and the right-hand side as
\begin{equation*}
b := \boldsymbol{W}^{-T}\hat{A}^TC_\eta^{-1}s.
\end{equation*}
Equation (\ref{eq:fewha}) is solved by applying a few iterations of a Jacobi preconditioned conjugate gradient (PCG) algorithm. The PCG is started with an initial guess obtained from the previous solution, which we refer to as warm restart.

Important parameters for measuring the performance of a solver are the dimensions of the matrices involved, which are listed in Table \ref{tab:matrixSize}. These values are defined by the number of layers $L$,  wavefront sensors $W$, subapertures $n_s$ and wavelet scales $J_\ell$.

\begin{table}[ht]
	\renewcommand{\arraystretch}{1.3}
	\centering
	\small
	\begin{tabular}{lrcl}
		Discrete Wavelet Transform $W^{-1}$ & $2^{2J_\ell}L$ & \hspace{-0.2cm}$\times$\hspace{-0.2cm} & $2^{2J_\ell}L$ \\
		SH Operator $\Gamma$ & $2n_s^2W$ & \hspace{-0.2cm}$\times$\hspace{-0.2cm} & $(n_s+1)^2 W$\\
		Bilinear Interpolation $P$ & $(n_s+1)^2 W$ & \hspace{-0.2cm}$\times$\hspace{-0.2cm} & $2^{2J_\ell}L$\\
		Left-hand side operator $M$ & $2^{2J_\ell}L$ & \hspace{-0.2cm}$\times$\hspace{-0.2cm} & $2^{2J_\ell}L$\\
	\end{tabular}
	\caption{Dimensions of the FEWHA matrices.}
	\label{tab:matrixSize}
\end{table}

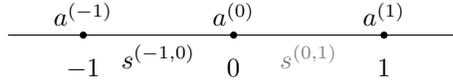
\begin{figure}
	\centering
	\begin{tikzpicture}[scale=2.0]
	\draw (-1.5,0)-- (1.5,0);
	\draw [fill=black] (-1,0) circle (0.02);
	\draw [fill=black] (0,0) circle (0.02);
	\draw [fill=black] (1,0) circle (0.02);
	\draw (0,0) node[above] {$a^{(0)}$};
	\draw (-1,0) node[above] {$a^{(-1)}$};
	\draw (1,0) node[above] {$a^{(1)}$};
	\draw (0.5,0) node[below, gray] {$s^{(0,1)}$};
	\draw (-0.5,0) node[below] {$s^{(-1,0)}$};
	\draw (-1,-0.1) node[below] {$-1$};
	\draw (0,-0.1) node[below] {$0$};
	\draw (1,-0.1) node[below] {$1$};
	\end{tikzpicture}
	\caption{Two-step delay}
	\label{fig:twostepdelay}
\end{figure}

There is a delay in the AO system between the time at which measurements are obtained and the time at which the correction is applied. Since the atmosphere is changing rapidly, a control scheme has to predict the DM shape update based on the collected measurements and the previous DM shape(s) \cite{Poettinger_2019}. We denote by the superscript indices $(-1)$, $(0)$ and $(1)$ the previous, the current and the next step of the loop. The superscript indices $(-1,0)$ and $(0,1)$ denote the measurements between the respective time steps.  Within FEWHA a two-step delay is used \cite{Yu14}, meaning that the new mirror shapes $a^{(1)}$ are determined from the reconstruction based on the measurements $s^{(-1,0)}$, which we simply denote by $s$, and from the previous mirror shapes $a^{(0)}$. The measurements $s^{(0,1)}$ are not available at time Step $1$. See Figure \ref{fig:twostepdelay} for a graphical representation. 

\begin{algorithm}[ht]
	\caption{FEWHA reconstruction algorithm}
	\small
	\label{alg:FEWHA}
	\begin{algorithmic}
		\STATE{\textbf{Input:}\quad\quad~$s=(s_g)^G_{g=1}$ (measurement vector)\\
			\quad\quad\quad\quad\quad~$gain$ (scalar weight)\\
			\quad\quad\quad\quad\quad~$c^{(0)}$ (previous wavelet coefficients)\\
			\quad\quad\quad\quad\quad~$b^{(0)}$ (previous right-hand side)\\
			\quad\quad\quad\quad\quad~$r^{(0)}$ (previous residual vector)\\
			\quad\quad\quad\quad\quad~$a^{(-1)}, a^{(0)}$ (previous two DM shape)}
		\STATE{\textbf{Output:}\quad~$a^{(1)}$ (next DM shape)}
		
		\IF{loop = closed}
		\STATE $s = s + \mathbf \Gamma a^{(-1)}$\label{alg.LTAO_MOAO.applyH}
		\ENDIF
		\vspace{0.3cm}
		
		\STATE $b^{(1)} = \mathbf W^{-T} {\mathbf A}^T \mathbf C_{\eta}^{-1} s$ \label{alg.LTAO_MOAO.RHS}
		\STATE $\bar{r} = b^{(1)} - \mathbf M c^{(0)} = (b^{(1)} - b^{(0)}) + r^{(0)}$ \label{alg.LTAO_MOAO.res} 	
		\vspace{0.3cm}				
		\STATE$(c^{(1)}, r^{(1)}) = \mathrm{PCG}(c^{(0)}, \bar{r})$ \label{alg.LTAO_MOAO.PCG} 
		\vspace{0.3cm}
		\STATE$\tilde a = \mathbf F \mathbf W^{-1} c^{(1)}$			\label{alg.LTAO_MOAO.P}
		
		\vspace{0.3cm}
		\IF{loop = closed}
		\STATE $a^{(1)} = a^{(0)} + \mathrm{gain} \cdot (\tilde a - a^{(-1)})$ \label{alg.LTAO_MOAO.control_LTAO}
		\ELSIF{loop = open}
		\STATE $a^{(1)} = (1 - \mathrm{gain}) \cdot a^{(0)} + \mathrm{gain} \cdot \tilde a$ \label{alg.LTAO_MOAO.control_LTAO_open}
		\ENDIF
	\end{algorithmic}
\end{algorithm}

The wavelet reconstructor is outlined in Algorithm \ref{alg:FEWHA}. The main input of FEWHA is the measurement vector $s$ and the output is the new shape of the mirrors $a^{(1)}$. An AO system can either operate in closed loop or in open loop. If the AO system is running in open loop, the measurements are obtained directly from the wavefronts. If closed loop control is applied, the DMs correct the wavefront before the measurements are obtained. In this case the pseudo open loop measurements are computed in a first step (see line \ref{alg.LTAO_MOAO.applyH}) of the algorithm as the sum of the actual residual measurement (stemming from the WFS) and the simulated SH measurements through the DM(s), i.e., the measurements correspond to the residuals of the corrected wavefront after the DM correction. Due to the two-step delay, the DM shape from the previous step is used. In line \ref{alg.LTAO_MOAO.RHS} the right-hand side (RHS) $b^{(1)}$ is computed with the new measurement vector and, subsequently, the residual vector $r^{(0,1)}$ is updated in line \ref{alg.LTAO_MOAO.res}. The atmospheric reconstruction takes place in line \ref{alg.LTAO_MOAO.PCG}, i.e., the PCG algorithm is applied to Equation~\eqref{eq:fewha}. In line \ref{alg.LTAO_MOAO.P} the layers are fitted to actuator commands. In lines \ref{alg.LTAO_MOAO.control_LTAO} and \ref{alg.LTAO_MOAO.control_LTAO_open} the closed or open loop control is applied. The new DM shapes are calculated as the linear combination of the current and the reconstructed DM shapes. The scalar weight between those two quantities is called gain and has a value between zero and one. With such a gain control the stability of the reconstruction is improved. For closed loop control the artificially added DM shapes are subtracted from the computed mirror shapes, such that the $(a-a^{(-1)})$ corresponds to the reconstruction from the closed loop measurements.

\section{Real-time implementation of FEWHA}\label{sec:realtimecomp}
This section is dedicated to the optimization techniques we use to implement FEWHA on real-time hardware. From a theoretical point of view, we follow the work in \cite{Parallel_Computing} for CPUs and \cite{Cuda} for GPUs. Some optimizations techniques listed here were already shown within \cite{Yu14}. We continue this work and provide an optimized parallel implementation of FEWHA on CPUs and a completely new version of the algorithm for GPUs.

\subsection{Hardware architecture} 	
In general, three basic technologies are used for the real-time control of large telescopes: CPUs (Central Processing Unit), GPUs (Graphics Processing Unit), and FPGAs (Field Programmable Gate Array). In this work, we address the CPU as well as the NVIDIA GPU technology for the real-time implementation of FEWHA. Due to its high development costs the FPGA technology is omitted and planned as a future task.

A GPU is optimized for a rapid processing of simple tasks. Due to their highly parallel structure GPUs outperform CPUs by orders of magnitude in algorithms that can be easily split in numerous identical computational tasks to be processed in parallel. GPUs have a much higher memory bandwidth compared to CPUs, because more transistors are dedicated to data processing rather than caching or flow control. In contrast, for the CPU memory access latencies are avoided through large data caches and flow control. Generally, CPUs are managed by the operating system (OS). These OS often generates unintended side effects in latency, jitter and determinism of the control behavior for CPUs as well as GPUs \cite{Rodr15_ao4elt4_real-time-control}. 

\subsection{Matrix-free implementation}\label{sec:matrix-free}
Almost all components of FEWHA have an efficient matrix-free representation, which allows a significant reduction in storage and floating point operations. The most important and time consuming operators for FEWHA are the wavelet transform, the bilinear interpolation and the SH operator. All these operators have a sparse structure and the matrix-free representation leads to a linear computational complexity and memory requirement. Details about the matrix-free representation of the operators can be found in \cite{Yu14}. In \cite{StaAO4ELT} details about the number of floating point operations and memory usage for FEWHA applied to an ELT-sized problem are listed. 

\subsection{Parallelization}\label{sec:ParProgMo}
Without parallelization it would not be possible to meet the real-time requirements of a large AO system as, e.g., required for the ELT. FEWHA allows two types of parallelization, which we refer to as global and local parallelization. By global parallelization we understand the block decomposition of the algorithm as a whole, whereas local parallelization refers to the parallelization of the individual operators. The combination of these two strategies leads to a very efficient parallelization scheme, which we implement using dynamic parallelism \cite{DynamicParallel}. 

For parallelization on CPU we use the popular environment OpenMP, which is a portable standard for shared memory programming. The OpenMP API provides compiler directives, environmental variables as well as library routines for the programming languages C, C++ and Fortran. OpenMP is supported by several compilers as, e.g., GCC or Clang. An OpenMP program starts with the execution of a single thread, which runs the program sequentially until a parallel construct is reached. The code inside the parallel region is then executed concurrently by all the threads in the team. On the NVIDIA GPU we use CUDA to parallelize FEWHA. CUDA offers special C++ functions, called kernels, that are executed in parallel by $N$ threads. Each thread has a unique index that is accessible within the kernel. The thread index is a three dimensional vector, to provide a natural way of applying operations on vectors, matrices and tensors in parallel often needed in scientific computations. For recent GPUs there is a limit of $1024$ threads per block, because all threads are located on the same core and must share the memory resources there. A kernel can be executed by several blocks consisting up to $1024$ threads. Blocks are organized in a one-, two- or three-dimensional grid. The architecture of NVIDIA GPUs is built on an array of Streaming Multiprocessors (SMs). When a kernel is invoked by a CUDA program, the blocks of the grid are distributed to multiprocessors with available execution capacity. All threads of a thread block execute simultaneously on one multiprocessor, and multiple blocks can execute concurrently on one multiprocessor.

\subsubsection{Global parallelization}\label{sec:globalParallelization}
\begin{figure}[ht]
	\centering
	\includegraphics[width=0.9\textwidth]{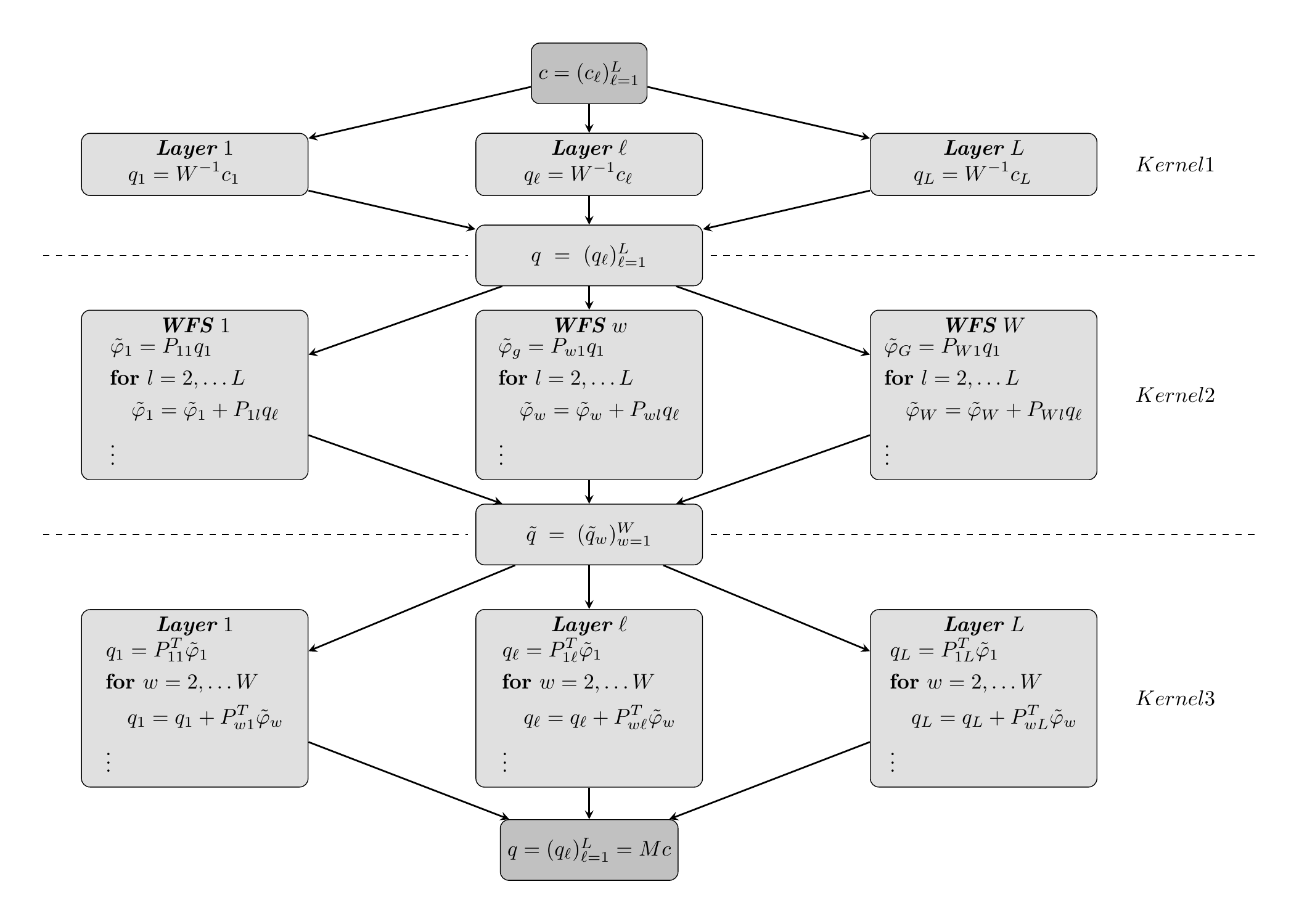}
	\caption{Parallelization of operator $M$.}
	\label{fig:parallelM}
\end{figure}

We illustrate the main idea behind the global parallelization of FEWHA on the left-hand side operator $M$ of Equation \eqref{eq:fewha} in Figure \ref{fig:parallelM}. This operator is applied once per PCG iteration and is the computationally heaviest part of the algorithm. Similar strategies are used for the right-hand side of Equation \eqref{eq:fewha}. Parallelization is applied over $L$ layers and $W$ WFS. $Kernel 1$ applies the inverse wavelet transform and is parallelizable over the number of layers $L$. $Kernel 2$ applies the atmospheric tomography operator $A$, which is composed into a SH-matrix $\Gamma$ and a bilinear interpolation matrix $P$, the inverse covariance matrix $C_\eta^{-1}$ and the transposed SH-matrix $\Gamma^T$ in parallel over the number of WFS $W$. Finally, $Kernel 3$ applies the transposed bilinear interpolation matrix $P^T$, the transposed inverse wavelet transform and adds the regularization term $\alpha D$ in parallel over the number of layers $L$. The algorithm is not perfectly parallelizable and after a certain number of steps, i.e., after each kernel, synchronization is necessary (highlighted with dashed lines). This is not optimal for certain hardware architectures, specially for GPUs, since they are optimized for computational throughput and latency is a problem.

On the CPU we define one parallel OpenMP region with $N := \max(L, W)$ number of threads for the program part that is responsible for the application of $M$ and run each loop inside in parallel. The reason for only one parallel region with $N$ threads instead of several regions with $L$ or $W$ threads is that creating threads is time consuming and empirical test show a better performance with this approach.

In a first attempt, we utilized unified memory for the parallel version on the GPU, because it provides an easy way to port an existing C++ code to CUDA. However, tests show that this approach is not satisfying and we decided to handle the copy operations from CPU to GPU and vise versa by our own. In fact, for the final version of FEWHA all computations are done by the GPU to minimize host/device memory transfer. We utilize CUDA kernels with $L$ or $W$ blocks, where each block has a single thread, for the parallel execution of the corresponding loops. The number of threads and blocks were optimized within several tests.

\subsubsection{Local parallelization}
The matrix-free operators (Section \ref{sec:matrix-free}) of FEWHA share the same basic structure, namely they consist of operations that are applied to a grid of values and can be computed independently from each other. This is exactly what Single Instruction Multiple Data (SIMD) is about. The main idea behind SIMD is to apply the same arithmetic operation to several data elements simultaneously using proper vector instructions. Such computations are supported by many CPUs as well as GPUs. The benefit of these computations is that a single vector instruction handles the computation of several elements, whereas a scalar instruction treats only one data element. Special vector load instructions are used to load the vector registers with the data from main memory. There are two common ways of dealing with vectorization: the use of vectorizing compilers and the use of a programming language. If there is no dependence from one loop iteration to the following, vectorizing compilers can transform them into an equivalent vector statement. We use explicit vectorization to parallelize the discrete wavelet transform, the SH operator and the bilinear interpolation of FEWHA. 

For the CPU we implement the local parallelization of the operators using Intel Thread Building Blocks (TBB) \cite{ProTBB}. TBB is a widespread C++ library for shared memory parallel programming. It provides broad support for parallelism that is beyond the scope of the C++ standard, although, with C++ 17 some of the features were included \cite{C++17}. We also tried an approach with OpenMP and dynamic parallelism, however, compared to SIMD and even compared to the non-parallelized version this shows a poor performance. 

On the GPU we implement the discrete wavelet transform as described, e.g., in \cite{DWT1} or \cite{DWT2}. Both the bilinear interpolation and the SH operator operate on a grid of values, where the computation on each value is independent from the others. Thus, they can be easily parallelized using CUDA dynamic parallelism with an optimized number of blocks and threads. Moreover, we apply a variety of optimizations techniques that are available on NVIDIA GPUs. Global memory loads and stores of the 32 threads of a warp are merged by NVIDIA GPUs into the fewest possible number of transactions. This effect is known as memory coalescing. We utilize this feature and align the data for the operators to minimize the DRAM bandwidth usage. In addition, CUDA provides user defined vectorized memory load and store instructions. Using them it was possible to improve the performance of the operators even further. These operations load and store data in 64- or 128-bit widths, and thus reduce the total number of instructions and latency, and improve bandwidth utilization. The fact that global memory and shared memory have different speed in access by several orders of magnitude offers another possibility to improve performance. We utilize shared memory for the SH operator to avoid slow global memory access. 

\begin{figure}[ht]
	\centering
	\includegraphics[width=1.0\textwidth]{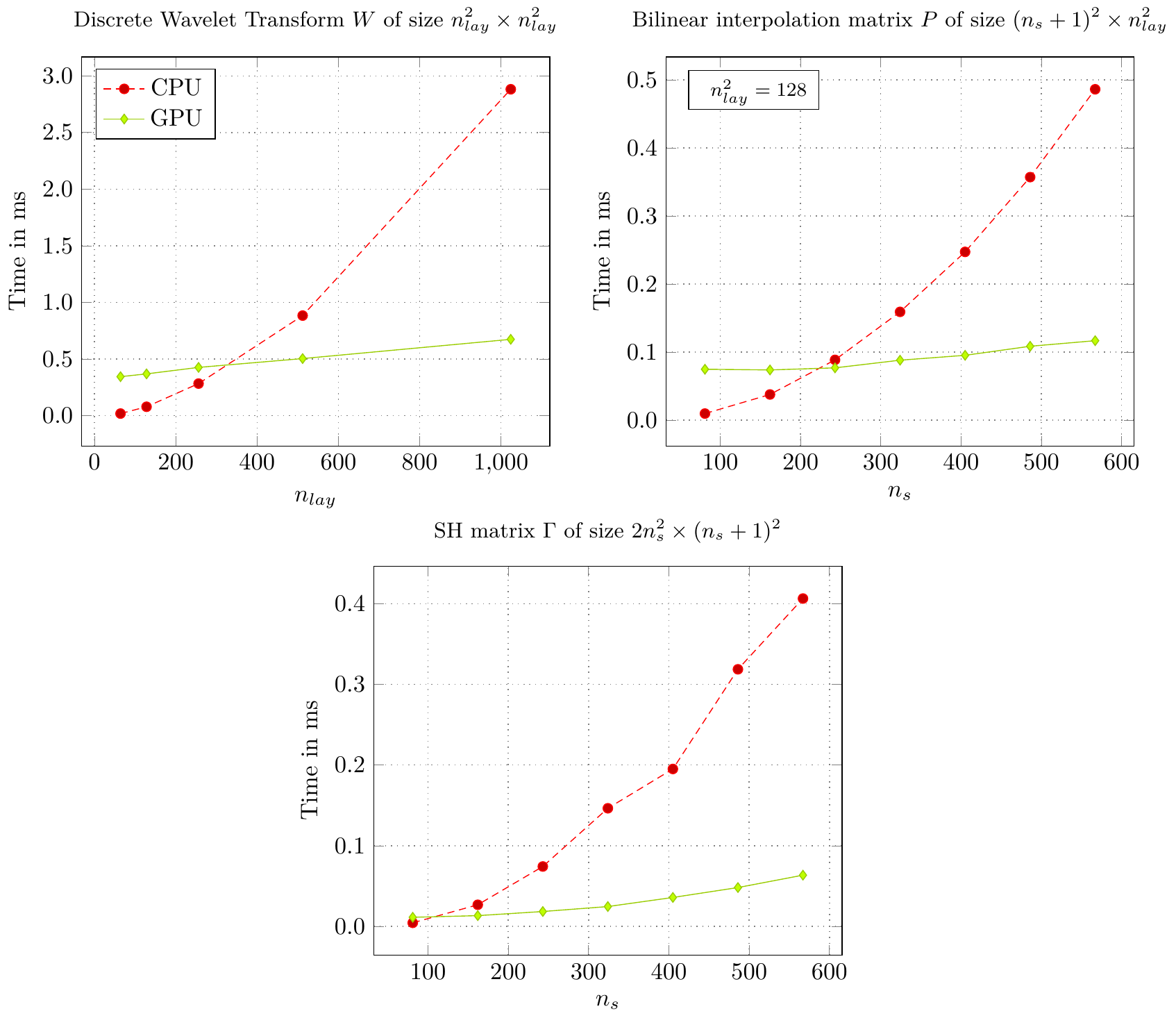}
	\caption{Timing of FEWHA operators.}
	\label{fig:operatorsTiming}
\end{figure}

Figure \ref{fig:operatorsTiming} shows the performance of the final, optimized, matrix-free versions of the operators on the CPU (red) and on the GPU (green). The computations are executed on the high-performance cluster Radon 1 for the CPU or a Tesla V100 GPU (for details about the hardware specifications we refer to Section \ref{sec:hardware}).  The relevant parameters for performance evaluation are the number of layer discretization points $n_{lay}$ and the number of subapertures $n_s$, which both are related to the dimension of the matrices. The upper left graph shows the performance of the discrete wavelet transform $W$ with matrix size $n_{lay}^2 \times n_{lay}^2$. On the upper right side the performance of the discretized bilinear interpolation operator $P$ with dimension $(n_s+1)^2 \times n_{lay}^2$ is shown and in the bottom we can see the performance of the discretized SH operator $\Gamma$. Although all these operators are applied in a matrix-free manner, the dimensions of the matrices still influence the computational speed. We observe that for small matrix sizes the CPU outperforms the GPU. However, for bigger matrix sizes the GPU clearly shows its benefits. The time needed for the GPU versions grows almost linearly with the matrix sizes. This is because of the very good parallelization of GPUs and their enormous computational throughput. For the bilinear interpolation $n_{lay}^2$ is fixed to $128$, because this value does not influence the overall performance as the number of entries per row is constant. Moreover, we can examine the benefit of shared memory at the graph of the SH operator. It is the only case where the GPU version for small dimensions is almost as fast as the CPU based implementation. 

\subsection{Preconditioned Conjugate Gradient Method on GPU}
Based on the work in \cite{PCG_GPU} we use an improved version of the PCG on the GPU as illustrated in Algorithm \ref{alg:PCG}. Our aim is to minimize the number of kernel calls, global memory loads and the communication overhead. By $(\cdot,\cdot)$ we denote the standard $\ell_2$ scalar product. These scalar products are grouped together and computed within a single kernel to reduce the number of synchronizing steps and kernels on GPU and CPU (steps 4-5). Within the main loop we load vectors at the same place and apply vector operations within a single kernel to allow multiple operations to reuse the data (steps 9-12). In the end, this approach has more floating point operations than the original one, however, shows a better performance on the GPU.  Within FEWHA a modified Jacobi preconditioner is utilized with a different weighting of the low and high frequency regimes (see \cite{YuHeRa13} for details). Apart from reducing the number of iterations, the benefit of using such an approach is an increased robustness and stability of the overall method. Moreover, this preconditioner can be formulated as a diagonal matrix, which is very efficient regarding the computational costs.

\begin{algorithm}[ht]
	\caption{PCG on GPU for $Mc=b$}
	\small
	\label{alg:PCG}
	\begin{algorithmic}
		\STATE{\textbf{Input:}\quad\quad~$r$ (residual vector)\\
			\quad\quad\quad\quad\quad~$J^{(-1)}$ (Jacobi preconditioner)\\
			\quad\quad\quad\quad\quad~$M$   (FEWHA left-hand side operator)\\
			\quad\quad\quad\quad\quad~$c$   (wavelet coefficient vector)\\
			\quad\quad\quad\quad\quad~$maxIter$ (max. number of PCG iterations)\\
			\quad\quad\quad\quad\quad~$a^{(-1)}, a^{(0)}$ (previous two DM shape)}
		\STATE{\textbf{Output:}\quad~$c$ (new wavelet coefficient vector)\\
			\quad\quad\quad\quad\quad~$r$ (new residual vector)}
		
		\FOR{$iter = 0,1,..., maxIter$} 
		\STATE $z = J^{-1}r$
		\STATE $s = Mz$
		\vspace{0.2cm}
		\STATE $\rho = (r,z)$
		\STATE $\mu = (s,z)$
		\vspace{0.2cm}
		\STATE $\beta = \rho / \rho_{old}$
		\STATE $\alpha = \rho / (\mu-\rho\beta/\alpha)$
		\STATE $\rho_{old} = \rho$
		\vspace{0.2cm}
		\STATE $p = z + \beta p$
		\STATE $q = s + \beta q$
		\STATE $c = c + \alpha p$
		\STATE $r = r - \alpha q$
		\ENDFOR
	\end{algorithmic}
\end{algorithm}
\newpage
\section{Numerical results}\label{sec:numerics}
\subsection{AO system configuration}
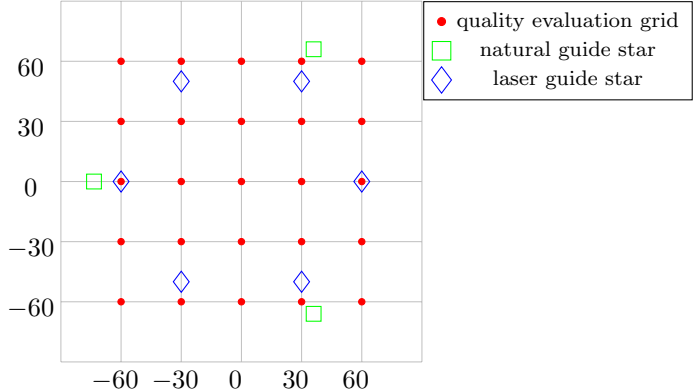
\begin{figure}
	\centering
	\begin{tikzpicture}[scale=0.8]
	\draw [very thin, lightgray] (0,0) grid (6,6);
	\draw (-0.5,1.2) node[below] {$-60$};
	\draw (-0.5,2.2) node[below] {$-30$};
	\draw (-0.5,3.2) node[below] {$0$};
	\draw (-0.5,4.2) node[below] {$30$};
	\draw (-0.5,5.2) node[below] {$60$};
	\draw (0.9, 0) node[below] {$-60$};
	\draw (1.9, 0) node[below] {$-30$};
	\draw (2.9, 0) node[below] {$0$};
	\draw (3.9,0) node[below] {$30$};
	\draw (4.9,0) node[below] {$60$};
	\foreach \x in {1,2,3,4,5}
	\draw [red] plot [only marks, mark size=1.5, mark=*] coordinates {(\x,1) (\x,2) (\x,3) (\x,4) (\x,5)};
	\draw [blue] plot [only marks, mark size=5, mark=diamond] coordinates {(2,4.666) (4,4.666) (1, 3) (5, 3) (2,1.333) (4,1.333)};
	\draw [green] plot [only marks, mark size=3.5, mark=square] coordinates {(0.55,3) (4.2,0.8) (4.2,5.2)};
	
	\begin{customlegend}[
	legend entries={
		quality evaluation grid,
		natural guide star, 
		laser guide star,
	},
	legend style={at={(10.5,6)},font=\footnotesize}]
	\addlegendimage{only marks, mark=*, color=red, mark size=1.5}
	\addlegendimage{only marks, mark=square, color=green, mark size=3.5}
	\addlegendimage{only marks, mark=diamond, color=blue, mark size=5}
	\end{customlegend}
	
	\end{tikzpicture}
	\caption{Star asterism of NGS, LGS and the quality evaluation grid.}
	\label{fig:MCAO}
\end{figure}

For all simulations we use a test configuration similar to MAORY \cite{MAORY}, which is an adaptive optics module for the ELT operating in MCAO. The telescope diameter is $39~m$, of which about $28~\%$ are obstructed. A SH-WFS is assigned to each guide star. We use six high resolution sensors for full wavefront correction with $80 \times 80$ subapertures and three low resolution sensors for tip-tilt aberrations (two with $1 \times 1$ and one with $2 \times 2$ subapertures). The LGS and NGS are located at a separation angle of $2$ and $8/3$ arcmin diameter, respectively. The quality for MCAO is measured in $25$ probe-star directions. The LGS and NGS star asterism with the quality evaluation grid for MCAO are illustrated in Figure~\ref{fig:MCAO}. The DMs are modeled by a piecewise bilinear function, with the grid nodes co-located with the actuator locations. For our analysis of the MCAO system we assume that the number of layers $L$ is equal to the number of DMs $M$. In this case, no interpolation onto the mirrors is required and the fitting operator is simply the identity operator, which improves the computational performance. If we reconstruct more layers than DMs ($L > M$) the quality increases, however, we have to solve an additional minimization problem, which is costly in terms of speed. See \cite{Yu14} for more details about mirror fitting. The method is configured to reconstruct up to $9$ layers of the simulated atmosphere using up to $9$ DMs. We vary the number of PCG iterations between $4$ (speed set up) and $8$ (quality set up). The dimension of this problem is very large, thus, a simple matrix-vector based approach for solving the corresponding atmospheric tomography problem is not feasible and iterative solvers are preferred. Table \ref{tab:maory_val} summarizes all important parameters for evaluating the computational performance. 

\begin{table}[ht]\label{tab:parameterSetting}
	\centering\renewcommand{\arraystretch}{1.3}\small
	\begin{tabular}{|l|l|l|}
		\hline
		\textbf{Description} & \textbf{Variable} & \textbf{Value} \\\hline\hline
		Number of wavefront sensors &  $W$ & $9$\\
		\quad Number of subapertures for WFS 1-6 & $n_{s1}^2$ & $80^2$\\
		\quad Number of subapertures for WFS 7 & $n_{s2}^2$ & $2^2$\\
		\quad Number of subapertures for WFS 8-9 & $n_{s3}^2$ & $1^2$\\
		Number of LGS & $G_{LGS}$ & $6$\\
		Number of NGS & $G_{NGS}$ & $3$\\
		Number of layers & $L$ & $3-9$\\
		Number of layer discretization points & $n_{lay}^2$ & $128^2$\\
		Number of DMs & $M$ & $3-9$\\
		\quad Number of actuators for DM 1 & $n_{a1}^2$ & $81^2$\\
		\quad Number of actuators for DM 2-5 & $n_{a2}^2$ & $48^2$\\
		\quad Number of actuators for DM 6-9 & $n_{a3}^2$ & $54^2$\\
		Number of PCG iterations & $n_{iter}$ & $4-8$ \\\hline
	\end{tabular}
	\caption{Parameter setting for performance evaluation.}
	\label{tab:maory_val}
\end{table}

\subsection{Hardware configuration}\label{sec:hardware}
We run the parallel CPU implementation of FEWHA on Radon1, the high performance cluster of the Radon Institute for Computational and Applied Mathematics in Linz. For our numerical simulations we use one compute node of Radon 1, that has two 8-core Intel Haswell processors (Xeon E5-2630v3, 2.4Ghz) and 128 GB of memory. 

For the performance tests of the GPU based FEWHA implementation we use a Tesla V100 GPU with CUDA 10.1. The Tesla V100 is a high-end GPU from NVIDIA optimized for deep learning and high-performance computing. Listing \ref{lst:deviceQuery} shows the output of the CUDA device query sample, provided by the CUDA toolkit. This example enumerates the properties of the available CUDA devices in the system. Important parameters for performance considerations are the clock rate, the maximal number of blocks per multiprocessor, the number of CUDA cores and the maximal number of threads per multiprocessor.

\begin{lstlisting}[float, basicstyle=\tiny\ttfamily, frame=single, caption={CUDA device query output.},captionpos=b, label={lst:deviceQuery}]
CUDA Device Query (Runtime API) version (CUDART static linking)

Detected 1 CUDA Capable device(s)

Device 0:"Tesla V100-PCIE-32GB"
CUDA Driver Version / Runtime Version          10.1 / 10.1
CUDA Capability Major/Minor version number:    7.0
Total amount of global memory:                 32480 MBytes (34058272768 bytes)
(80) Multiprocessors, ( 64) CUDA Cores/MP:     5120 CUDA Cores
GPU Max Clock rate:                            1380 MHz (1.38 GHz)
Memory Clock rate:                             877 Mhz
Memory Bus Width:                              4096-bit
L2 Cache Size:                                 6291456 bytes
Maximum Texture Dimension Size (x,y,z)         1D=(131072), 2D=(131072, 65536), 3D=(16384, 16384, 16384)
Maximum Layered 1D Texture Size, (num) layers  1D=(32768), 2048 layers
Maximum Layered 2D Texture Size, (num) layers  2D=(32768, 32768), 2048 layers
Total amount of constant memory:               65536 bytes
Total amount of shared memory per block:       49152 bytes
Total number of registers available per block: 65536
Warp size:                                     32
Maximum number of threads per multiprocessor:  2048
Maximum number of threads per block:           1024
Max dimension size of a thread block (x,y,z): (1024, 1024, 64)
Max dimension size of a grid size    (x,y,z): (2147483647, 65535, 65535)
Maximum memory pitch:                          2147483647 bytes
Texture alignment:                             512 bytes
Concurrent copy and kernel execution:          Yes with 7 copy engine(s)
Run time limit on kernels:                     No
Integrated GPU sharing Host Memory:            No
Support host page-locked memory mapping:       Yes
Alignment requirement for Surfaces:            Yes
Device has ECC support:                        Enabled
Device supports Unified Addressing (UVA):      Yes
Device supports Compute Preemption:            Yes
Supports Cooperative Kernel Launch:            Yes
Supports MultiDevice Co-op Kernel Launch:      Yes
Device PCI Domain ID / Bus ID / location ID:   0 / 0 / 9
\end{lstlisting}

\subsection{Computational performance}
\begin{figure}[ht]
	\centering
	\captionsetup{justification=centering}
	\begin{minipage}{.49\textwidth}
		\centering
		\includegraphics[width=0.9\textwidth]{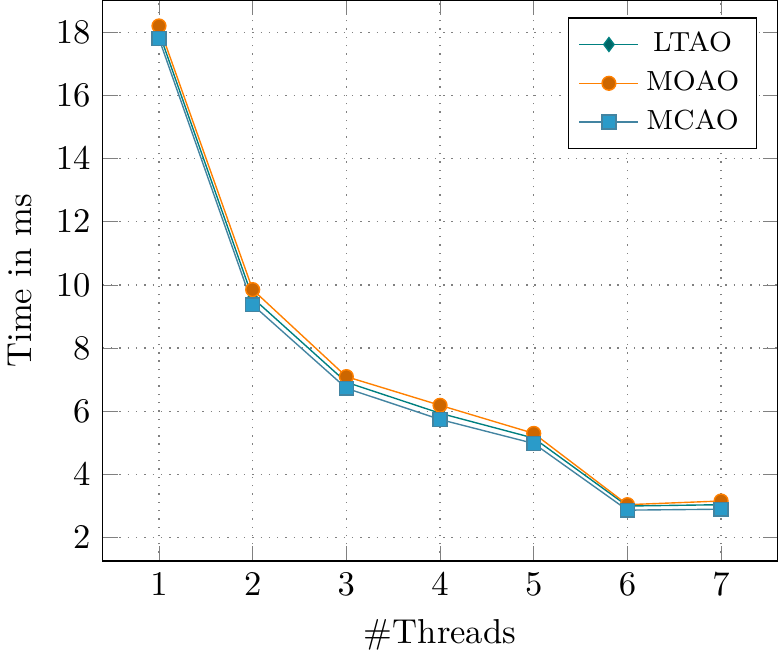}
		\captionof{figure}{Scalability of FEWHA with\\different number of threads on CPU for $L=M=6$ and $4$ PCG iterations.}
		\label{fig:FEWHAtimingThreads}
	\end{minipage}
	\begin{minipage}{.49\textwidth}
		\centering
		\includegraphics[width=0.9\textwidth]{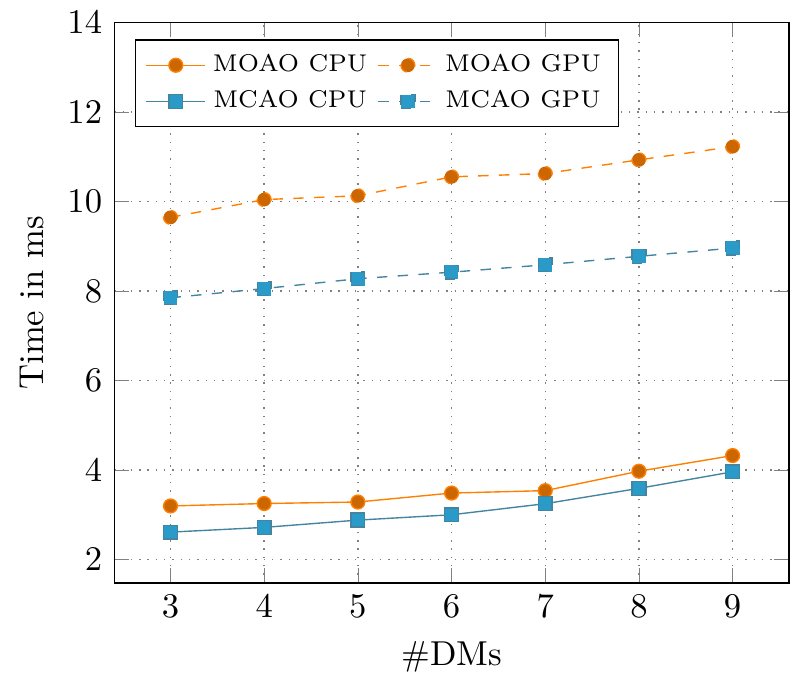}
		\captionof{figure}{Scalability of FEWHA with\\different number of DMs (equal to number of layers) and $4$ PCG iterations.}
		\label{fig:FEWHAtimingDM}
	\end{minipage}
\end{figure}

In this section we demonstrate the computational performance of FEWHA on the hardware described in the previous section for the test setting defined in Table \ref{tab:maory_val}. For our first test runs we fix the number of layers and DMs to $6$ and the number of PCG iterations to $4$. The parallel optimized CPU version takes about $3$ ms for each reconstruction whereas the optimized GPU implementation needs $8.6$ ms. On the CPU the algorithm performs best with only $6$ threads, which is related to the number of layers, as shown in Figure \ref{fig:FEWHAtimingThreads}. For a higher number of threads the elapsed time stays the same or even increases. We obtain a better performance if the number of threads is set for each parallel region separately, to either $6$ for parallelization over the number of layers or $9$ for parallelization over the number of WFS. The poor performance of the GPU is not totally surprising. Looking at Figure \ref{fig:operatorsTiming} we can observe that, using the MAORY parameter settings, we are on the lower end of the graphs for which the CPU is a little faster than the GPU. Moreover, for the global parallelization scheme we are just parallelizing over $6$ layers, or $9$ wavefront sensors, which is far too little to fully utilize the GPU resources. The work which is then executed in each parallel block (related to the boxes in Figure \ref{fig:parallelM}) is quite complicated and not a simple operation applied to a big amount of data. To obtain best performance on the CPU we only need $6$ or $9$ threads, which indicates that the level of parallelism is not very high. 

\begin{figure}[ht]
	\centering
	\captionsetup{justification=centering}
	\includegraphics[width=0.9\textwidth]{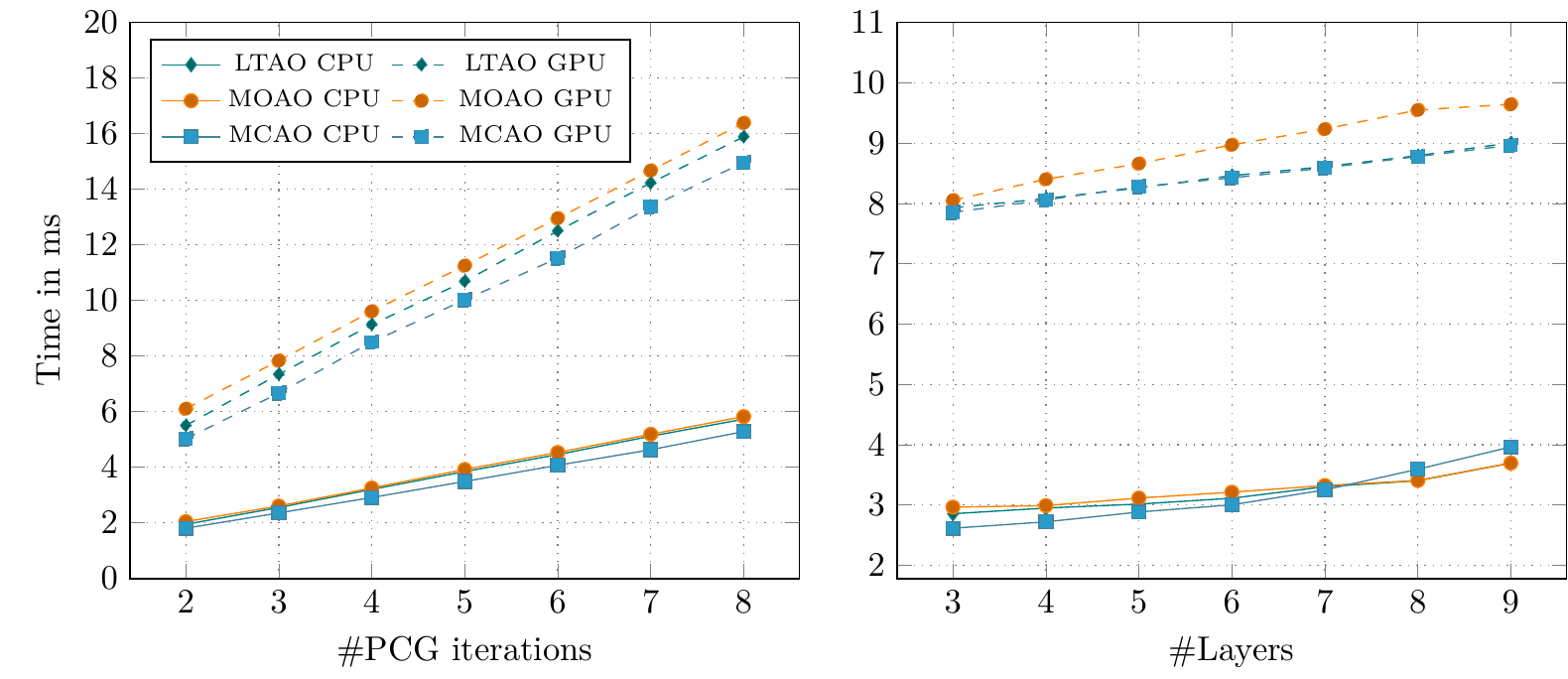}
	\caption{Scalability of FEWHA with PCG iterations (left) and layers (right). We simulated the left graph with $L=M=6$ and the right one with $4$ PCG iterations.}
	\label{fig:FEWHAtimingPCG}
\end{figure}

\begin{figure}[ht]
	\centering
	\captionsetup{justification=centering}
	\includegraphics[width=0.8\textwidth]{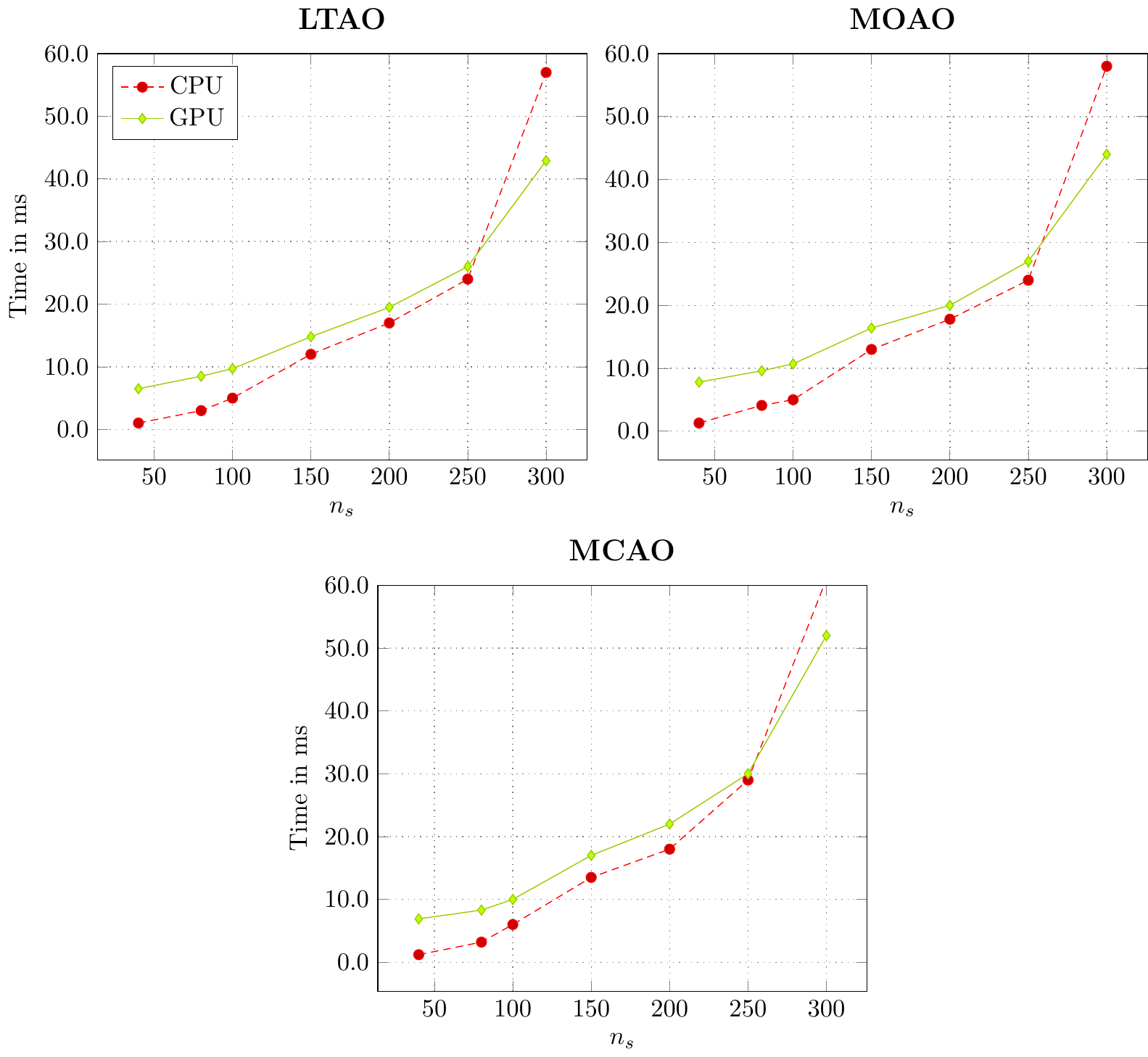}
	\caption{Scalability of FEWHA for different number of subapertures with $4$ PCG iterations and $L=M=6$.}
	\label{fig:FEWHAtimingSubap}
\end{figure}

Based on these investigations we start to analyze the behavior of the algorithm for the three different AO systems LTAO, MOAO and MCAO. Moreover, we vary the PCG iterations, number of layers and DMs. In the LTAO tests we use the same star asterism and WFS as for MCAO but with a single DM. Figure \ref{fig:FEWHAtimingPCG} shows the behavior for different number of PCG iterations (on the left side) and a different number of layers (on the right side). In Figure \ref{fig:FEWHAtimingDM} we vary the number of DMs for the MOAO and MCAO operating modes, since LTAO operates with a single DM. All test cases show a similar performance, namely that the CPU clearly outperforms the GPU. The low level of global parallelism influences the performance of the GPU for all our settings. Even if we go beyond $9$ layers or DMs, we are still at a very low level of global parallelism. The PCG method is not parallelizable, thus, having more iterations is clearly less effective. Looking closer at Figure \ref{fig:operatorsTiming} we examine that the local parallelization strategy on the GPU is more efficient for a higher number of subapertures for all the operators involved. We observe a similar behavior for FEWHA, as illustrated in Figure \ref{fig:FEWHAtimingSubap}. We increase the dimension of the setting up to $n_s = 300$ subapertures, and define the number of actuators for DM 1 as $n_{act} = n_s + 1$ and the number of layer discretization points $n_{lay}$ as the next higher power of two. For more than $250$ subapertures the GPU starts to show its benefits and beats the CPU in terms of speed for all three operating systems.  

\begin{figure}[ht]
	\centering
	\captionsetup{justification=centering}
	\begin{minipage}{.49\textwidth}
		\centering
		\includegraphics[width=1.\textwidth]{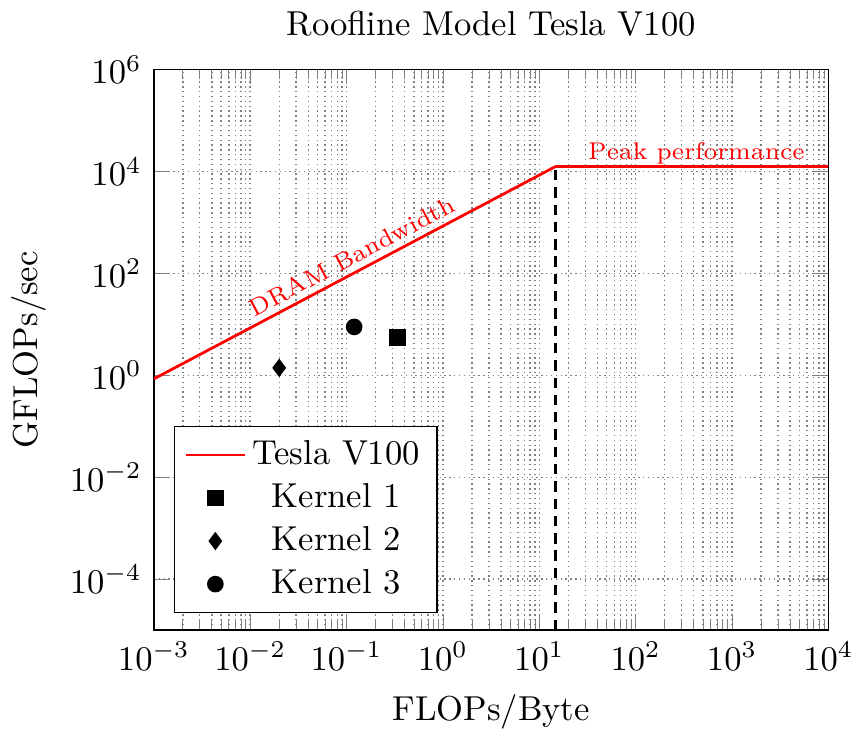}
	\end{minipage}
	\begin{minipage}{.49\textwidth}
		\centering
		\includegraphics[width=1.\textwidth]{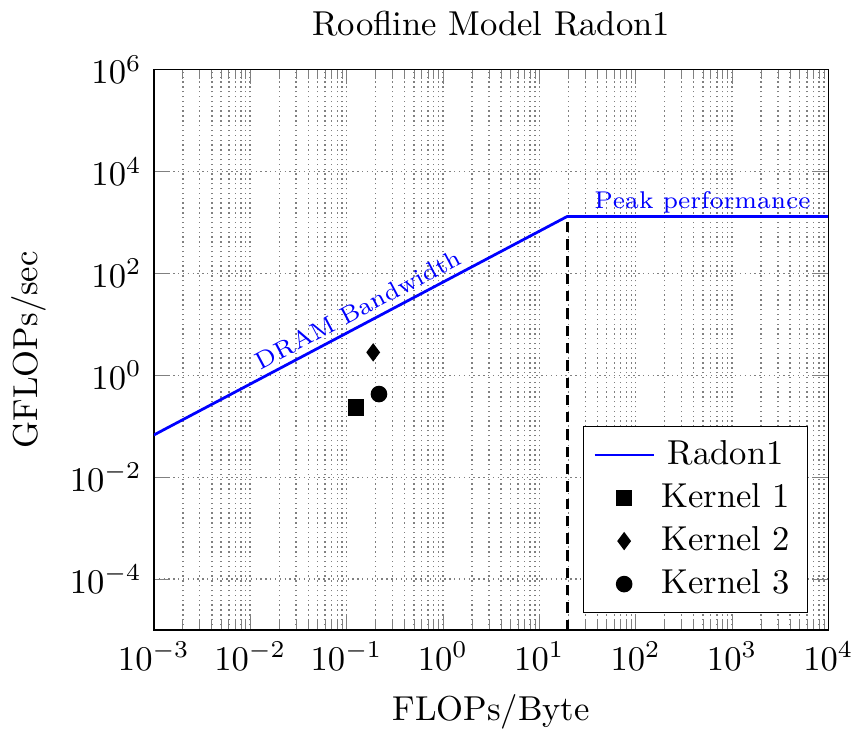}
	\end{minipage}
	\caption{Roofline model for the NVIDIA Tesla V100 GPU (on the left) and one node of Radon1 cluster (on the right) with the main kernels from the global parallelization scheme of the operator $M$.}
	\label{fig:RooflineModels}
\end{figure}

All these observations suggest that our application is bandwidth bounded, i.e., we can not optimize the performance by further parallelization or vectorization as memory accesses are the bottleneck. The roofline model provides a nice way to visualize the trade-off between computational intensity and data movement \cite{williams2009}. Hence, it offers us the possibility to identify if the memory bandwidth is really the limiting factor for our algorithm. The aim of this model is to show the computational performance, memory bandwidth and memory locality in one chart for the kernels of interest. Figure \ref{fig:RooflineModels} shows the roofline model for the Tesla V100 (on the left) and one node of the Radon1 cluster (on the right) for the main kernels of the global parallelization strategy of operator $M$ for the MAORY test case with $4$ PCG iterations. Note, that the outcome of the CPU and GPU kernels is the same, however, the implementation differs and is optimized for the given hardware, thus, the dots within the two plots are not at the same place. For a detailed description of the kernels see Section \ref{sec:globalParallelization}. For the GPU we use NVIDIA Nsight Compute to determine the measurements for the roofline model and for the CPU we use the Intel Advisor. Both profiling tools offer an intuitive way to create roofline models either for NVIDIA or Intel hardware. The x-axis shows the arithmetic intensity measured in FLOPs/Byte and the y-axis the performance in GFLOPs/sec, both in logarithmic intensity. The diagonal line represents how many bytes of data a given memory can deliver per second, whereas the horizontal line represents the number of floating point computations per second the given hardware can perform. Kernels that lie on left side of the dashed, black line are memory bandwidth bounded, whereas kernels that lie on the right side are compute bounded. We can observe that all our three kernels lie within the memory bandwidth bounded area.

All these test runs are based on a matrix-free implementation. For a matrix-based version of FEWHA we would have to store a huge amount of data, which leads to time consuming copy operations between CPU and GPU for a MAORY-sized test problem. Moreover, a (sparse) matrix representation is even more memory bandwidth bounded.

\section{Conclusion and Future Work}
In this paper, we continued the work of \cite{Yu14, HeYu13,YuHeRa13, YuHeRa13b} on the iterative solver FEWHA by considering a parallel implementation. We assessed the performance of this method and studied how well it can be parallelized on different computational platforms. As expected, the choice of the computational architecture has a crucial impact on the performance of the algorithm. There are various applications in different scientific areas that can gain a significant speed-up from GPUs. However, we demonstrated by numerical simulations that for ELT-sized problems FEWHA performs better on CPUs. Moreover, the development effort is much higher for GPUs. It is easy to learn CUDA, but it requires experience to program efficient code which fully utilizes GPU resources. The mathematical methodologies used within FEWHA allow to solve the atmospheric tomography problem with a very low number of floating point operations. This is a huge benefit on CPUs, however, GPUs are made for computational throughput, thus, are not the optimal architecture to solve an ELT-sized problem with FEWHA. Nevertheless, we showed that for an increasing number of subapertures and actuators, the GPU tends to outperform CPUs. In future work we aim to assess the performance of FEWHA on FPGAs. Additionally, we plan to analyze the feasibility of reduction methods for a matrix-based version of FEWHA, i.e., without the iterative PCG method and with a matrix-vector multiplication approach instead. 

\section{Acknowledgments}
The project has received funding from the European Union’s Horizon 2020 research and innovation programme under the Marie Sk\l odowska-Curie Grant Agreement No. 765374.
 
\FloatBarrier
\bibliography{eso_sorted}
\bibliographystyle{abbrv}
\end{document}